\documentclass[reqno]{amsart}

\usepackage{amsmath, amssymb}
\usepackage[mathscr]{eucal}
\usepackage{hyperref}
\hypersetup{colorlinks=true, linkcolor=blue, urlcolor=blue, citecolor=[rgb]{0, 0.8, 0}, linktocpage=true}

\theoremstyle{plain}
\newtheorem*{theorem}{Theorem}

\numberwithin{equation}{section}

\newcommand{\du}{\mathrm{d}}

\DeclareMathOperator{\ind}{ind}

\title[A Riesz basis criterion for Schr\"{o}dinger operators]{A Riesz basis criterion for Schr\"{o}dinger operators with boundary conditions dependent on the eigenvalue parameter}

\author{Namig J. Guliyev}
\address{Institute of Mathematics and Mechanics, Azerbaijan National Academy of Sciences, 9 B.~Vahabzadeh str., AZ1141, Baku, Azerbaijan.}
\email{njguliyev@gmail.com}



\begin{document}
\maketitle
\begin{abstract}
We establish a criterion for a set of eigenfunctions of the one-dimensional Schr\"{o}dinger operator with distributional potentials and boundary conditions containing the eigenvalue parameter to be a Riesz basis for $\mathscr{L}_2(0,\pi)$.
\end{abstract}

\setcounter{tocdepth}{1}
\tableofcontents

\section{Introduction and main result} \label{sec:introduction}

In this paper we continue the study of one-dimensional Schr\"{o}dinger operators with distributional potentials and boundary conditions containing rational Herglotz--Nevanlinna functions of the eigenvalue parameter initiated in~\cite{G19a}. These operators are generated by the differential equation
\begin{equation} \label{eq:SL}
  - \left( y^{[1]} \right)'(x) - s(x) y^{[1]}(x) - s^2(x) y(x) = \lambda y(x)
\end{equation}
and the boundary conditions
\begin{equation} \label{eq:boundary}
  \frac{y^{[1]}(0)}{y(0)} = -f(\lambda), \qquad \frac{y^{[1]}(\pi)}{y(\pi)} = F(\lambda),
\end{equation}
where $s \in \mathscr{L}_2(0, \pi)$ is real-valued, $y^{[1]}(x) := y'(x) - s(x) y(x)$ denotes the \emph{quasi-derivative} of $y$ with respect to $s$, and
\begin{equation} \label{eq:f_F}
  f(\lambda) = h_0 \lambda + h + \sum_{k=1}^d \frac{\delta_k}{h_k - \lambda}, \qquad F(\lambda) = H_0 \lambda + H + \sum_{k=1}^D \frac{\Delta_k}{H_k - \lambda}
\end{equation}
are rational Herglotz--Nevanlinna functions, i.e., $h_0, H_0 \ge 0$, $h, H \in \mathbb{R}$, $\delta_k, \Delta_k > 0$, $h_1 < \ldots < h_d$, $H_1 < \ldots < H_D$. Our aim in this paper is to prove a criterion for (a subset of) the eigenfunctions of this boundary value problem to be a Riesz basis for $\mathscr{L}_2(0,\pi)$.

In~\cite{G17}, to each function $f$ of the form~(\ref{eq:f_F}) we assigned its \emph{index} (an integer) which, roughly speaking, counts the number of poles of this function. More precisely, each finite pole is counted twice and a pole at infinity (if any) once:
\begin{equation*}
  \ind f := \begin{cases} 2 d + 1, & h_0 > 0, \\ 2 d, & h_0 = 0. \end{cases}
\end{equation*}
This notion allowed us in that paper (see also~\cite{G18a} and \cite{G19a}) to formulate various direct and inverse spectral results for boundary value problems with boundary conditions of the form~(\ref{eq:boundary}), (\ref{eq:f_F}) in a unified manner, i.e. without considering separate cases as it is usually done in the literature. Define a nonnegative integer $N$ by
\begin{equation*}
  N := \left\lceil \frac{\ind f}{2} \right\rceil + \left\lceil \frac{\ind F}{2} \right\rceil,
\end{equation*}
where the ceiling function $\lceil \cdot \rceil$ denotes the smallest integer not smaller than the argument. Let $\Theta := \{ n_1, \ldots, n_N \} \subset \mathbb{N} \cup \{ 0 \}$ be a set of $N$ distinct indices. Then from the asymptotics of the eigenvalues and the eigenfunctions one sees that the sequence $\{ \psi_n \}_{n \notin \Theta}$ of appropriately normalized (see below) eigenfunctions asymptotically behaves as (and is, in fact, quadratically close to) the orthonormal basis
\begin{equation*}
  \sqrt{\frac{2}{\pi}} \cos \left( \left( n + N - \frac{\ind f + \ind F}{2} \right) x + \frac{\ind f}{2} \pi \right)
\end{equation*}
(which coincides with one of $\{ \cos n x \}_{n=0}^{\infty}$, $\{ \cos (n + 1/2) x \}_{n=0}^{\infty}$, $\{ \sin (n + 1/2) x \}_{n=0}^{\infty}$, or $\{ \sin (n + 1) x \}_{n=0}^{\infty}$ up to a constant factor). Hence it seems reasonable to hope that the sequence $\{ \psi_n \}_{n \notin \Theta}$ will be a Riesz basis for $\mathscr{L}_2(0,\pi)$, i.e. the image of an orthonormal basis under a bounded invertible operator~\cite{B51}. This is indeed the case when the boundary conditions do not contain the eigenvalue parameter at all or only one of them depends on the eigenvalue parameter, and can easily be established by using the transformation operators~\cite[Theorem 6.2]{HM04}, \cite[Corollary 5.1]{G05}. However, in the general case it is quite possible for $\{ \psi_n \}_{n \notin \Theta}$ not to be a Riesz basis for $\mathscr{L}_2(0,\pi)$. It turns out that whether this sequence is a (Riesz) basis or not depends on the invertibility of a certain $N \times N$ matrix defined in terms of the spectral characteristics of the boundary value problem~(\ref{eq:SL})-(\ref{eq:boundary}).

We need some additional definitions to state our result. Denote by $W$ the diagonal matrix with diagonal entries $\delta_1^{-1}$, $\ldots$, $\delta_d^{-1}$, $h_0^{-1}$, $\Delta_1^{-1}$, $\ldots$, $\Delta_D^{-1}$, $H_0^{-1}$, where the $(d+1)$-th entry (respectively, the last entry) is omitted whenever $h_0 = 0$ (respectively, $H_0 = 0$), and consider the Hilbert space $\mathcal{H} = \mathscr{L}_2(0,\pi) \oplus \mathbb{C}^N$ with inner product given by
\begin{equation*}
  \left\langle \begin{pmatrix} y \\ \widehat{y} \end{pmatrix}, \begin{pmatrix} z \\ \widehat{z} \end{pmatrix} \right\rangle_{\mathcal{H}} := \int_0^{\pi} y(x) \overline{z(x)} \,\du x + \widehat{z}^{\dagger} W \widehat{y},
\end{equation*}
where the superscript $^{\dagger}$ denotes the conjugate transpose. Most of our matrices will have real entries and for them the conjugate transpose coincides with the ordinary transpose. The boundary value problem~(\ref{eq:SL})-(\ref{eq:boundary}) is equivalent to an eigenvalue problem for a self-adjoint operator in $\mathcal{H}$ with discrete spectrum (see~\cite[Section 2.3]{G19a} for details), in the sense that they both have the same eigenvalues $\lambda_n$ and this operator has an orthonormal basis of eigenvectors of the form
\begin{equation} \label{eq:psi_vector}
  \begin{pmatrix} \psi_n \\ \widehat{\psi}_n \end{pmatrix},
\end{equation}
where $\psi_n$ are eigenfunctions of~(\ref{eq:SL})-(\ref{eq:boundary}) and
\begin{equation*}
  \widehat{\psi}_n := \begin{pmatrix} \frac{\delta_1 \psi_n(0)}{\lambda_n - h_1} & \dots & \frac{\delta_{d} \psi_n(0)}{\lambda_n - h_{d}} & -h_0 \psi_n(0) & \frac{\Delta_1 \psi_n(\pi)}{H_1 - \lambda_n} & \dots & \frac{\Delta_{D} \psi_n(\pi)}{H_{D} - \lambda_n} & H_0 \psi_n(\pi) \end{pmatrix}^{\dagger}
\end{equation*}
(with the obvious modifications when one or both of $h_0$ and $H_0$ equal zero). We can (and will) choose $\psi_n$ to be real-valued. This kind of linearization procedure goes back at least to a 1956 book by Friedman \cite[pp. 205--207]{F56} and can even be generalized to arbitrary (not necessarily rational) Herglotz--Nevanlinna functions~\cite{G19b}. We define $M_{\Theta}$ as the matrix whose rows consist of (the entries of) the vectors $\widehat{\psi}_{n_k}$:
\begin{equation*}
  M_{\Theta} := \sum_{k=1}^N e_k \widehat{\psi}_{n_k}^{\dagger},
\end{equation*}
where $\{ e_k \}_{k=1}^N$ is the standard basis of $\mathbb{C}^N$. Our main result can now be stated as follows.

\begin{theorem} \label{thm:identities}
  The sequence $\{ \psi_n \}_{n \notin \Theta}$ is a Riesz basis for $\mathscr{L}_2(0,\pi)$ if and only if the matrix $M_{\Theta}$ is invertible.
\end{theorem}

We will prove this theorem in the next section. As remarked by Gelfand~\cite{G51}, to prove that a sequence is a Riesz basis it suffices to construct a new inner product equivalent to the original one, with respect to which this sequence becomes an orthonormal basis~\cite[Theorem 1.9]{Y01}. The main idea of our proof is to demonstrate that this new inner product in $\mathscr{L}_2(0,\pi)$ can be constructed in terms of the inner product of the space $\mathcal{H}$ in a straightforward way.

Now that we have this general result, the following question naturally arises: for a given problem, roughly speaking, what part of $N$-tuples $\Theta$ satisfies the condition of the theorem? Since, intuitively speaking, a generic matrix is invertible, one might expect that the share of $N$-tuples with $\det M_{\Theta} = 0$ will be negligible in some sense. Indeed, as we have already pointed out, if only one of the boundary conditions depends on the eigenvalue parameter then each $M_{\Theta}$ is invertible. On the other hand, however, for symmetric boundary value problems with linear dependence on the eigenvalue parameter ($s(x) + s(\pi-x) = 0$, $f = F$, and $1 \le \ind f \le 2$), roughly speaking, only half of all $N$-tuples satisfies the condition of the theorem. We will discuss these issues in Section~\ref{sec:cases}.

\section{Proof} \label{sec:proof}

We start with the ``only if'' part. If the matrix $M_{\Theta}$ is not invertible then the vectors $\widehat{\psi}_{n_k}$ are linearly dependent, i.e.,
\begin{equation*}
  \sum_{k=1}^N \alpha_k \widehat{\psi}_{n_k} = 0
\end{equation*}
for some $\alpha_k$, not all of them being zero. The function
\begin{equation*}
  y(x) := \sum_{k=1}^N \alpha_k \psi_{n_k}(x)
\end{equation*}
cannot be identically equal to zero, since otherwise the orthonormal system
\begin{equation*}
  \begin{pmatrix} \psi_{n_k} \\ \widehat{\psi}_{n_k} \end{pmatrix}, \qquad k = 1, \ldots, N
\end{equation*}
would also be linearly dependent. Moreover,
\begin{equation*}
  \int_0^{\pi} y(x) \psi_n(x) \,\du x = \left\langle \begin{pmatrix} y \\ 0 \end{pmatrix}, \begin{pmatrix} \psi_n \\ \widehat{\psi}_n \end{pmatrix} \right\rangle_{\mathcal{H}} = \sum_{k=1}^N \alpha_k \left\langle \begin{pmatrix} \psi_{n_k} \\ \widehat{\psi}_{n_k} \end{pmatrix}, \begin{pmatrix} \psi_n \\ \widehat{\psi}_n \end{pmatrix} \right\rangle_{\mathcal{H}} = 0
\end{equation*}
for every $n \notin \Theta$. Hence $y \neq 0$ is orthogonal to all the functions of the sequence $\{ \psi_n \}_{n \notin \Theta}$, and thus this sequence cannot be complete in $\mathscr{L}_2(0,\pi)$.

We now turn to the ``if'' part. Our immediate aim is to define a new inner product equivalent to the usual one in $\mathscr{L}_2(0,\pi)$ and such that the sequence $\{ \psi_n \}_{n \notin \Theta}$ is an orthonormal basis with respect to this new inner product. Since we already have the Hilbert space $\mathcal{H}$ and the orthonormal sequence~(\ref{eq:psi_vector}) in this space, the most straightforward way to achieve this is to map $\mathscr{L}_2(0,\pi)$ into $\mathcal{H}$ in such a way that $\{ \psi_n \}_{n \notin \Theta}$ are mapped to their corresponding vectors~(\ref{eq:psi_vector}) and then ``transfer'' the inner product on $\mathcal{H}$ to $\mathscr{L}_2(0,\pi)$. With this goal in mind, we define the mapping $y \mapsto y_{\Theta}$, $\mathscr{L}_2(0,\pi) \to \mathbb{C}^N$ by the formula
\begin{equation} \label{eq:y_Theta}
  y_{\Theta} := - W^{-1} M_{\Theta}^{-1} \sum_{k=1}^N e_k \int_0^{\pi} y(x) \psi_{n_k}(x) \,\du x.
\end{equation}
One can easily verify that $\left( \psi_n \right)_{\Theta} = \widehat{\psi}_n$ for $n \notin \Theta$. Now we introduce a new inner product in $\mathscr{L}_2(0,\pi)$ by the obvious expression
\begin{equation*}
  \langle y, z \rangle_{\Theta} := \left\langle \begin{pmatrix} y \\ y_{\Theta} \end{pmatrix}, \begin{pmatrix} z \\ z_{\Theta} \end{pmatrix} \right\rangle_{\mathcal{H}}.
\end{equation*}
It is trivial to check that this is indeed an inner product and $\{ \psi_n \}_{n \notin \Theta}$ are orthonormal with respect to it. That this inner product is equivalent to the original one follows from the inequalities
\begin{equation*}
  \int_0^{\pi} |y(x)|^2 \,\du x \le \langle y, y \rangle_{\Theta} \le \left( 1 + \| W^{-1} \| \| M_{\Theta}^{-1} \|^2 \sum_{k=1}^N \int_0^{\pi} |\psi_{n_k}(x)|^2 \,\du x \right) \int_0^{\pi} |y(x)|^2 \,\du x,
\end{equation*}
where $\| \cdot \|$ denotes the operator norm. It remains to verify the completeness. To this end, suppose that $\langle y, \psi_n \rangle_{\Theta} = 0$ for all $n \notin \Theta$. Then
\begin{equation*}
  \begin{pmatrix} y \\ y_{\Theta} \end{pmatrix} = \sum_{k=1}^N \alpha_k \begin{pmatrix} \psi_{n_k} \\ \widehat{\psi}_{n_k} \end{pmatrix}
\end{equation*}
for some $\alpha_k$, this being a consequence of the orthogonality of the vector on the left-hand side to the vectors from~(\ref{eq:psi_vector}) with $n \notin \Theta$. Therefore~(\ref{eq:y_Theta}) yields
\begin{equation*}
\begin{split}
  \sum_{k=1}^N \alpha_k \widehat{\psi}_{n_k} = y_{\Theta} = \sum_{k=1}^N \alpha_k \left( \psi_{n_k} \right)_{\Theta} &= \sum_{k=1}^N \alpha_k \left( \widehat{\psi}_{n_k} - W^{-1} M_{\Theta}^{-1} e_k \right) \\
  &= \sum_{k=1}^N \alpha_k \widehat{\psi}_{n_k} - W^{-1} M_{\Theta}^{-1} \sum_{k=1}^N \alpha_k e_k.
\end{split}
\end{equation*}
Since $M_{\Theta}^{-1}$ and $W^{-1}$ are both obviously invertible, we obtain $\alpha_k = 0$ for $k = 1$, $\ldots$, $N$, and hence $y = 0$.

\section{Some special cases} \label{sec:cases}

In general, it appears to be rather difficult to characterize the $N$-tuples $\Theta$ for which $M_{\Theta}$ is invertible in terms of $\Theta$ itself. Two particular cases when this is possible have already been studied (for continuously differentiable $s$) in the literature (see below for references). We now discuss these two cases.

\subsection{Dependence on the eigenvalue parameter only in one boundary condition} \label{ss:one}

In this section we assume that one of the boundary coefficients, say $F$, is constant. As we have noted in the introduction, in this case by using the transformation operators~\cite{HM04}, one can deduce that $M_{\Theta}$ is invertible for every $\Theta$. We now want to derive this result as a corollary of our theorem.

In our case, the matrix $M_{\Theta}$ (after obvious cancellations) has either the form
\begin{equation} \label{eq:matrices}
  \begin{pmatrix}
    p_1 \left( \lambda_{n_1} \right) & \dots & p_d \left( \lambda_{n_1} \right) & p \left( \lambda_{n_1} \right) \\
    \vdots & \ddots & \vdots & \vdots \\
    p_1 \left( \lambda_{n_d} \right) & \dots & p_d \left( \lambda_{n_d} \right) & p \left( \lambda_{n_d} \right) \\
    p_1 \left( \lambda_{n_{d+1}} \right) & \dots & p_d \left( \lambda_{n_{d+1}} \right) & p \left( \lambda_{n_{d+1}} \right)
  \end{pmatrix} \qquad \text{or} \qquad \begin{pmatrix}
    p_1 \left( \lambda_{n_1} \right) & \dots & p_d \left( \lambda_{n_1} \right) \\
    \vdots & \ddots & \vdots \\
    p_1 \left( \lambda_{n_d} \right) & \dots & p_d \left( \lambda_{n_d} \right)
  \end{pmatrix},
\end{equation}
depending on whether $h_0 > 0$ or $h_0 = 0$, where we denoted
\begin{equation*}
  p(\lambda) := \prod_{k=1}^d (h_k - \lambda), \qquad p_m(\lambda) := \prod_{\substack{k = 1 \\ k \ne m}}^d (h_k - \lambda).
\end{equation*}
One possible way to prove the invertibility of these matrices, as was done in~\cite{KA06}, is to reduce this problem to the invertibility of a Cauchy matrix whose determinant has a closed-form expression (the Cauchy double alternant identity).

However, since we do not need the value of this determinant, but only need to verify the invertibility of the matrices in~(\ref{eq:matrices}), one can also proceed as follows. If
\begin{equation*}
  \sum_{m=1}^d \alpha_m p_m \left( \lambda_{n_k} \right) + \alpha p \left( \lambda_{n_k} \right) = 0, \qquad k = 1, \ldots, d+1,
\end{equation*}
then the polynomial $\alpha_1 p_1 (\lambda) + \ldots + \alpha_d p_d (\lambda) + \alpha p (\lambda)$ of degree at most $d$ has $d+1$ roots and hence it must be identically zero. From the obvious identities $p_m(h_k) = 0$ for $m \neq k$ we obtain $\alpha_1 = \ldots = \alpha_d = \alpha = 0$ and thus the columns of the first matrix in~(\ref{eq:matrices}) are linearly independent. This proves the invertibility of this matrix. The case $h_0 = 0$ is similar.

\subsection{Linear dependence on the eigenvalue parameter} \label{ss:linear}

The case when one or both of the boundary conditions depend in a linear fashion on the eigenvalue parameter (i.e., $\max \{ \ind f, \ind F \} \le 2$ in our notation) is the most extensively studied case in the literature. We mention only the very recent paper~\cite{K19} and refer the reader to the bibliography therein. Since we have already discussed the case when only one of the boundary conditions depends on the eigenvalue parameter, we now assume that both of them contain the eigenvalue parameter (i.e., $\min \{ \ind f, \ind F \} \ge 1$).

In order to completely characterize those pairs $\Theta = \{ n_1, n_2 \}$ for which $M_{\Theta}$ is invertible, we need some more definitions from~\cite{G17}, \cite{G19a}. We assign to each rational Herglotz--Nevanlinna function $f$ of the form~(\ref{eq:f_F}) two polynomials $f_\uparrow$ and $f_\downarrow$ by writing this function as
\begin{equation*}
  f(\lambda) = \frac{f_\uparrow(\lambda)}{f_\downarrow(\lambda)},
\end{equation*}
where
\begin{equation*}
  f_\downarrow(\lambda) := h'_0 \prod_{k=1}^d (h_k - \lambda), \qquad h'_0 := \begin{cases} 1 / h_0, & h_0 > 0, \\ 1, & h_0 = 0. \end{cases}
\end{equation*}
For each $n \in \mathbb{N} \cup \{ 0 \}$ we define $\beta_n \ne 0$ as the unique number for which
\begin{equation*}
  \chi_n(x) = \beta_n \varphi_n(x),
\end{equation*}
where $\varphi_n$ and $\chi_n$ are the (necessarily linearly dependent) eigenfunctions of~(\ref{eq:SL})-(\ref{eq:boundary}) satisfying the conditions $\varphi_n(0) = f_\downarrow(\lambda_n)$, $\varphi_n^{[1]}(0) = -f_\uparrow(\lambda_n)$, $\chi_n(\pi) = F_\downarrow(\lambda_n)$, and $\chi_n^{[1]}(\pi) = F_\uparrow(\lambda_n)$. Using the results of~\cite{G18a} and \cite{G19a} one can deduce that these numbers have alternating signs and asymptotically behave as
\begin{equation} \label{eq:beta}
  \beta_n = (-1)^n \left( n - \frac{\ind f + \ind F}{2} \right)^{\ind F - \ind f} \left( 1 + \xi_n \right), \qquad \xi_n \in \ell_2.
\end{equation}

We now return to the main topic of this subsection. Under our assumptions, the matrix $M_{\Theta}$ (again after obvious cancellations) has the form
\begin{equation*}
  \begin{pmatrix}
    1 & \beta_{n_1}^{-1} \\
    1 & \beta_{n_2}^{-1}
  \end{pmatrix}.
\end{equation*}
Thus the invertibility of $M_{\Theta}$ is equivalent to the condition $\beta_{n_1} \ne \beta_{n_2}$. On the basis of the discussion in the preceding paragraph, we immediately conclude that $M_{\Theta}$ is always invertible for indices $n_1$ and $n_2$ of different parity. For indices of same parity, both cases are possible. As one extreme example, we have $\beta_n = (-1)^n$ for all $n$ when the boundary value problem is symmetric (i.e., $s(x) + s(\pi-x) = 0$ and $f = F$) and hence $M_{\Theta}$ is never invertible for $n_1$ and $n_2$ of same parity. On the other hand, using the inverse spectral theory developed in~\cite{G19a}, one can produce a boundary value problem of the form~(\ref{eq:SL})-(\ref{eq:boundary}) for arbitrary sequence of distinct numbers $\beta_n$, as long as they satisfy the requirements of the preceding paragraph. In the case of the latter boundary value problem, $M_{\Theta}$ is invertible for all pairs of indices $n_1$ and $n_2$.

As a final observation, we note that it is also possible to obtain some results of asymptotic character when $\ind F \ne \ind f$. For example, given $n_1$, $M_{\Theta}$ is invertible for all sufficiently large $n_2$. In the case of summable potentials (i.e., for absolutely continuous $s$), one can say even more: $M_{\Theta}$ is invertible for all sufficiently large $n_1$ and $n_2$. In other words, there can be only finitely many pairs $\Theta = \{ n_1, n_2 \}$ for which $M_{\Theta}$ is not invertible, since in this case $n \xi_n \to 0$ in~(\ref{eq:beta}).

\end{document}